# A Little Reflection about the Sleeping Beauty Problem

André Luiz Barbosa [1]
http://www.andrebarbosa.eti.br – andrelb808@gmail.com
Non-commercial projects: Simple – PLC Simulator & LCE – Electric Commands Language

***Abstract***. *This paper presents a little reflection about the Sleeping Beauty Problem, maybe contributing to shed light on it and perhaps helping to find a simple and elegant solution that could definitively resolve the controversies about it.*

Mathematics Subject Classification (2020). Primary 91B06; Secondary 62C12.

Keywords. Sleeping Beauty Problem, Probability Theory, Self-Locating Belief, Philosophy of Science, Decision Theory, Cognitive Science.

## 1. Introduction

The Sleeping Beauty Problem [1,2], also known as the thorny issue of self-locating belief, is a thought experiment in Probability Theory, Philosophy of Science, and related fields. The problem can be stated as follows:

"Sleeping Beauty (we shall call her SB) volunteers to undergo the following experiment and is told all of the following details: On Sunday she will be put to sleep. Once or twice, during the experiment, SB will be awakened, interviewed, and put back to sleep with an amnesia-inducing drug that makes her forget that awakening. A fair coin will be tossed to determine which experimental procedure to undertake:

- If the coin comes up heads, SB will be awakened and interviewed on Monday only.
- If the coin comes up tails, she will be awakened and interviewed on Monday and

---

[1] Bachelor of Computer Science from UFG, in Brazil.

Tuesday.

In either case, she will be awakened on Wednesday without interview and the experiment ends. Any time SB is awakened and interviewed she will not be able to tell which day it is or whether she has been awakened before. During the interview SB is asked: "What is your credence (degree of belief) now for the proposition that the coin landed heads?"

The question is: What should SB's degree of belief (credence) be about the outcome of the coin toss when she is awakened?" [1]

There are two main competing viewpoints on this problem:

1. The "Thirders" Viewpoint [3,7]: Thirders argue that SB should assign a probability of **1/3** to the coin landed heads. They argue that she should consider all **3** possible awakenings with interview (that they think that are equally likely (equiprobable) events on each awakening, by the Principle of Indifference [5,6]), where if the coin landed heads, then will be only **1** awakening with interview (on Monday).

2. The "Halfers" Viewpoint [4,8]: Halfers argue that SB should assign a probability of **1/2** to the coin landed heads, regardless of whether she is awakened on Monday or Tuesday. They argue that she should reason as if she is flipping a fair coin on each awakening with interview and that her degree of belief (credence) should be updated accordingly (that is, when SB awakes, without new information, her degree of belief (credence) should be the same [**1/2**] on Monday and Tuesday as on Sunday). Notice that her degree of belief (credence) is not a mere arbitrary subjective concept, but rather rational, objective and mathematical probability [9].

Hence, we can see that the Sleeping Beauty Problem is quite controversial, and there is no consensus on the correct answer. It has sparked lively discussions in Philosophy, Probability Theory, Decision Theory and Cognitive Science, with different interpretations and arguments put forth by proponents of each point of view.

Nevertheless, here we do discuss a little reflection upon it, proposing two slightly different variants of the problem (that we shall call "*Double Sleeping Beauty Problem*" in Section 2, and "*Infinite Double Sleeping Beauty Problem*" in Section 3), and prove a theorem on the matter, maybe contributing to shed light on it, perhaps even helping the researchers in the area to find a simple and elegant solution that could definitively resolve those controversies. That is, if this short paper has any value, then it is hoped that it can help clarify that superb problem and contribute to a better understanding of it.

## 2. The Double Sleeping Beauty Problem

**Definition 2.1: The Double Sleeping Beauty Problem.** Two Sleeping Beauties (we shall call them SB1 and SB2) volunteer to undergo the following experiment and are told all of the following details: On Sunday they will be put to sleep in separate rooms. Once or twice, during the experiment, the Sleeping Beauties will be awakened, interviewed, and put back to sleep with an amnesia-inducing drug that makes their forget that awakening. A fair coin will be tossed to determine which experimental procedure to undertake:

- If the coin comes up heads, then SB1 will be awakened and interviewed on Monday only, whereas SB2 will be awakened and interviewed on Monday and Tuesday.

- If the coin comes up tails, then SB2 will be awakened and interviewed on Monday only, whereas SB1 will be awakened and interviewed on Monday and Tuesday.

In either case, they will be awakened on Wednesday without interview and the experiment ends.

Any time Sleeping Beauties are awakened and interviewed they will not be able to tell which day it is or whether they have been awakened before. During the interview of SB1 she is asked: "What is your credence (degree of belief) now for the proposition that the coin landed heads?"; whereas during the interview of SB2 she is asked: "What is your credence (degree of belief) now for the proposition that the coin landed tails?".

The question is: What should each SB1's and SB2's degree of belief (credence) be about the respective outcome of the coin toss when they are awakened?

It seems now, with this modified version of that experiment, that the thirder position at original version of the problem is maybe inconsistent: Following that answer that the degree of belief (credence) of the lonely SB should be **1/3**, then in Double Sleeping Beauty Problem above both answers should be **1/3**, too, of course, since those SB1 and SB2 situations are symmetric and exactly the same as that the lonely SB faces in the original problem. However, this answer implies that the added probabilities that the coin landed either heads or tails would be equal to just **2/3** (**1/3 + 1/3**), which seems inconsistent with the actual value (**1**), since that experiment involves a unique classical fair coin toss that must certainly land exclusively either heads or tails (for these events are equally likely, mutually exclusive, and exhaustive).

On the other hand, notice, however, that the halfer position at original version of the problem seems now plainly consistent: Following that answer that the degree of belief (credence) of the lonely SB should be **1/2**, then in Double Sleeping Beauty Problem above both answers should be **1/2**, too, of course, and then this answer implies that the added probabilities that the coin landed either heads or tails would be equal to **1** (**1/2 + 1/2**), which seems now utterly consistent with the actual value.

Apparently, the possible mistake of the thirder position was to consider that the correct measure of the SB's degree of belief (credence) about the proposition that the coin landed heads should be the quantity of awakenings with interview with respect to this outcome (just one for it against two for the coin landed tails, equivocally supposed equiprobable events on each awakening). However, this quantity seems not more than a simple kind of imaginary 'prize' or 'reward' associated with a particular outcome, conceivably obliterating our judgment about the related probabilities, but unable to really represent the true SB's degree of belief (credence) about the proposition that that outcome has occurred:

After she awaken, maybe what SB should consider is not her possible quantity of awakenings with interview, but whether she is in the awakenings-chain Monday-Tuesday (probability = **1/2**) or in the unique-awakening Monday (probability = **1/2**, too), in order to evaluate her degree of belief (credence) about the proposition that the coin landed heads (**1/2**). Notice that, in other possible versions of the original problem, that awakenings-chain can be of arbitrary length, even infinite one (Monday-Tuesday, Monday-Tuesday-Wednesday, Monday-Tuesday-Wednesday-Thursday, …), but, in any case, the probability that she is in whatever awakenings-chain shall be always equal to **1/2**, exactly the same probability that she is in the unique-awakening Monday.

## 3. Rejecting Brevity

A reviewer of the previous version of this paper said: "The brevity of the submission makes it difficult to fully assess the merits and novelty of the proposed definition". So, let's challenge the brevity and think about another interesting variation of that problem.

**Definition 3.1: The Infinite Double Sleeping Beauty Problem.** Suppose that in the creation of the Universe, a fair coin was tossed to determine which course it would undertake:

- If the coin came up heads, then all the women *in*[2] the Earth would be immortal and live infinite number of lives, successive and perpetually being born, living, dying, and reincarnating, where before each reincarnation they would suffer an amnesia-inducing process that would make them to forget all their previous livings. On the other hand, all the men in the Earth would be mortal and live only one life, being born, living, dying, and done.

- If the coin came up tails, then all the men in the Earth would be immortal and live infinite number of lives, successive and perpetually being born, living, dying, and reincarnating, where before each reincarnation they would suffer an amnesia-inducing process that would make them to forget all their previous livings. On the other hand, all the women in the Earth would be mortal and live only one life, being born, living, dying, and done.

Consider herein that all the women and men in the Earth are absolutely rational and firmly believe that those statements in Definition 3.1 above are true. So, during the living of a married woman (in a heterosexual marriage) in the Earth she is asked: "What is your credence (degree of belief) now for the proposition that you are immortal (living an infinite succession of lives in the Earth)?"; and, at this same moment, her husband is asked too: "What is your credence (degree of belief) now for the proposition that you are immortal (living an infinite succession of lives in the Earth)?".

The question is: What should each those woman's and man's degree of belief (credence) be about the respective immortality?

It seems again, now with another modified version of that problem, that the thirder position at original version is newly maybe inconsistent: Following exactly the same reasoning that entails that the degree of belief (credence) of the lonely SB should be **1/3** (that the coin in the original problem landed heads), then in the problem above both answers, from woman and man, should be **1** (100% certainty in the own immortality), of course, since both would live an infinite quantity of lives in the Earth in a case (immortality) against only one in the alternative one (no immortality). However, this answer implies that the added probabilities that the coin tossed in the creation of the Universe landed either heads or tails would be equal to **2** (**1** + **1**), which seems inconsistent with the actual value (**1**), since that fair coin toss involves a unique classical one too: one of them, either the woman or the man, even though absolutely rational person, would be necessarily absolutely wrong.

Moreover, if the thirder point of view was correct, then, in the scope of the problem above, there would be sadly nearly four billion people living in the Earth, even though extremely rational ones, into complete delusional credence (degree of belief) about the absolute certainty with respect to their own immortality.

[2] The Earth is our home, not a mere surface: So, we live **in** it, not **on** it.

## 4. Simplified Mathematical Analysis

Maybe the possible fundamental equivocated reasoning of the "Thirders" could be an incorrect application of the Principle of Indifference [5,6], considering the three possible SB's awakenings with interview in the original Sleep Beauty Problem (Monday-Heads [$H_1$], Monday-Tails [$T_1$] and Tuesday-Tails [$T_2$]) as equally likely (equiprobable) events in *the SB's point of view*, where they are really not so: The probability of getting Monday-Heads, $P(H_1)$ = **1/2**, whereas $P(T_1) = P(T_2) =$ **1/4**, since $P(T_1)$ and $P(T_2)$ are equally likely events in *the SB's point of view*, and their added probabilities sum up to **1/2**, for they occur only when that coin landed tails, and $P(T)$ = **1/2** = $P(H) = P(H_1)$.

That is, the probability of getting tails [$P(T)$ = **1/2**] is *equally distributed* between Monday and Tuesday awakenings [$P(T_1)$ and $P(T_2)$] in *the SB's point of view* (for, when she awakes, without new information, she concludes that it can be either Monday or Tuesday, there is no way she knows which day she has wakened), whereas the probability of getting heads [$P(H)$] is *associated* only with Monday awakening [$P(H_1)$], where *the SB's point of view* is individualized herein by means that amnesia-induced that makes her forget the previous awakening.

However, if on Monday awakening it is told SB that today is Monday, then $P(T_1)$ *collapses* into **1/2**, $P(T_2)$ *collapses* into **0**, and $P(H_1)$, $P(H)$ and $P(T)$ *continue* being equal to **1/2**, since that *distribution* of $P(T)$ between $P(T_1)$ and $P(T_2)$ would no longer make sense in this case. And if on occasional Tuesday awakening is told SB that today is Tuesday, then $P(T_2)$ and $P(T)$ *collapse* into **1**, and $P(T_1)$, $P(H_1)$ and $P(H)$ *collapse* into **0**, by similar reasoning. Furthermore, notice that $P(H) + P(T) =$ **1** in all these cases.

Now, in order to turn the reasoning and visualization of these probabilities clearer, let's see them into the table below (upon *the SB's point of view*) [3]:

| **SB Does Not Know Which Day She Has Wakened** | | |
|---|---|---|
| Whichever Awakening Day with Interview | $P(H) = P(H_1) =$ **1/2** | P = 1/2; H; M \| T; P = 1/4 = 1/4; T |
| | $P(T_1) = 1/4$ and $P(T_2) =$ **1/4** | |
| | $P(T) = P(T_1) + P(T_2) =$ **1/2** | |
| | $P(H) + P(T) =$ **1** | |
| **SB Does Know Which Day She Has Wakened** | | |
| Monday Awakening Day | $P(H) = P(H_1) =$ **1/2** | |
| | $P(T_1) = 1/2$ and $P(T_2) =$ **0** | |
| | $P(T) = P(T_1) + P(T_2) =$ **1/2** | |
| | $P(H) + P(T) =$ **1** | |
| Tuesday Awakening Day | $P(H) = P(H_1) =$ **0** | |
| | $P(T_1) = 0$ and $P(T_2) =$ **1** | |
| | $P(T) = P(T_1) + P(T_2) =$ **1** | |
| | $P(H) + P(T) =$ **1** | |

**Table 4.1** Distributions and Collapses of SB's Point of View Probabilities

Notice that, despite $P(T_1)$ and $P(T_2)$ [**1/4**] are half of $P(H_1)$ [**1/2**] when SB does not know which day she has wakened, the averages number of times these events happen (their frequencies) are equal in the original Sleeping Beauty Problem, because the number of awakenings with interview is the double when the coin lands tails compared with heads, which compensates that difference.

[3] The picture into the table was adapted from [1] and is distributed under the CC BY-SA 4.0 license.

## 5. The Origin of that Thirder Position Confusion

Another reviewer of the previous version of this paper said, explaining the origin of that thirder position confusion:

"... Probabilists have rules for describing random experiments. When the experiment is set up according to those rules, there can be a couple of ways of looking at it, but the assumptions involved and conclusions are clear. I think that your intuition about the problem is essentially right. ...

The essential point is that the thirder position stems from a misunderstanding of how to set up the analysis of the experiment as a proper probability model and then some confusion between what is a *probability* and what is a *ratio of expectations*.

The sample space ***S*** of a random experiment is the set of possible outcomes. When the random experiment is conducted, one of those outcomes occurs. In the SB experiment, if we are going to make ***S*** = {HMon, TMon, TTue}, we have to be thinking of TMon and TTue as alternatives, not something that can both occur. This can make sense if we are willing to think of SB's waking-up day as a random element. But if we do this, there is no justification for saying that the probability measure on the sample space has P{HMon} = {PTMon} = P{TTue} = **1/3**. In fact, if we randomly choose SB's waking up day when we have a choice, then the probability measure consistent with the description of the experiment has P{HMon} = ½, with P{TMon} = P{TTue} = ¼.

Specifically, with the sample space and probability measure defined as above, the conclusion that P{TMon} = P{HMon} is incorrect. The description of the experiment states that 'a fair coin will be flipped' (which, using the above sample space, can be expressed as P{HMon} = ½). However, the fact that the coin is fair isn't used at all in the thirder argument that justifies P{HMon} = P{TMon} = P{TTue} = **1/3**. To see why it needs to be, ask yourself whether you would be comfortable in saying P{HMon} = P{TMon} = P{TTue} = **1/3** if the probability that the coin came up heads was **999999/1000000** [or any **1 - Ɛ**, for **1 > Ɛ > 0**].

I'll now explain what the thirder position illustrated, for example, by the simulations performed by Derek Muller at https://www.youtube.com/watch?app=desktop&v=XeSu9fBJ2sI actually calculates.

The only source of randomness is the coin toss, which can come up **H** or **T**. So, the sample space is ***S*** = {H, T}. If the coin comes up heads, Derek puts a **1** in the Monday Heads column, if it comes up tails, he puts **1**s in both the Monday Tails and Tuesday Tails columns. He then calculates the proportion of all **1**s that appear in the Monday Heads column. This is equivalent to defining a random variable **W** such that W(H) = **1** and W(T) = **2**.

This random variable has probability mass function P(W=1) = P(W=2) = ½.

The expected value E(W) = ½ * **1** + ½ * **2** = **3/2**, which makes sense – the average number of wakeups is **3/2**.

The number of wakeups that are equal to **1** is captured in the random variable **Y** that has Y(H) = **1** and Y(T) = **0**, which has probability mass function: P(Y=1) = P(Y=0) = ½.

and expected value E(Y) = ½ * **1** + ½ * **0 =** ½.

We can then observe that E(Y)/E(W) = **1/3**.

The crucial observation is that E(Y)/E(W) is a ratio of expectations, not a probability. The 'thirder position' makes the mistake of confusing these two things."

Then, as in Table 4.1 above, in order to turn that reasoning and visualization of those probabilities and ratios of expectations [of the quantity of SB's hits] clearer, let's see them into the table below (upon *the SB's point of view* too):

| **Coin Landed** | **Probability (SB's Degree of Belief or Credence)** | **Qty of Awakenings w.Interview** | **SB's Choice (Guess or Bet)** | **Qty of Hits** | **Expectation of Qty of Hits** | **Ratio of Expectations** |
|---|---|---|---|---|---|---|
| Heads | 1/2 | 1 | **Heads** | 1 | 1*1/2+0*1/2 | (1/2)/(1/2+1) |
| Tails | 1/2 | 2 | | 0 | = 1/2 | = **1/3** |
| Heads | 1/2 | 1 | **Tails** | 0 | 0*1/2+2*1/2 | (1)/(1/2+1) |
| Tails | 1/2 | 2 | | 2 | = 1 | = **2/3** |

**Table 5.1** Difference Between Probability and Ratio of Expectations Demonstrating Thirder Error

Additionally, in order to turn that reasoning and visualization of those probabilities and ratios of expectations yet clearer, let's see an alternative of them into the table below, considering now that SB is rewarded with a prize of **1** USD if that coin landed heads and she hits it (that is, betting that the coin landed heads), and with a prize of **2** USD if the coin landed tails and she hits it (that is, betting that the coin landed tails), leading to the fact that she should always bet that the coin landed tails, in order to wisely maximize her expectation of prize (expected prize), but not because she rationally believes that the probability (herein, SB's degree of belief or credence) of the coin landed tails is greater than it landed heads (since they are really equal, and SB knows that):

| **Coin Landed** | **Probability (SB's Degree of Belief or Credence)** | **Prize for Hit (in US$ or USD)** | **SB's Choice (Guess or Bet)** | **Prize** | **Expectation of Prize** | **Ratio of Expectations** |
|---|---|---|---|---|---|---|
| Heads | 1/2 | $ 1 | **Heads** | $ 1 | 1*1/2+0*1/2 | (1/2)/(1/2+1) |
| Tails | 1/2 | $ 2 | | $ 0 | = $ 1/2 | = **1/3** |
| Heads | 1/2 | $ 1 | **Tails** | $ 0 | 0*1/2+2*1/2 | (1)/(1/2+1) |
| Tails | 1/2 | $ 2 | | $ 2 | = $ 1 | = **2/3** |

**Table 5.2** Another Example of the Difference Between Probability and Ratio of Expectations

## 6. Sleeping Beauty's Theorem

Based on the above, we can even formalize and generalize the knowledge developed herein into a theorem on the matter:

**Theorem 6.1. Sleeping Beauty's Theorem.** *Suppose that the original Sleeping Beauty Problem is generalized in order to permit arbitrary positive [integer] numbers **n** and **m** of awakenings with interview depending on the outcome of the coin flipping: if it was heads, then there shall be **n** ones, otherwise (if it was tails), **m** ones. Then, even in this generalized* *case, SB should assign a probability of* ***1/2*** *to the coin landed heads (and, consequently,* ***1/2*** *to the coin landed tails, too), but now the ratios of expectations of the quantity of her hits shall be* ***n/(n+m)*** *when the SB's choice (guess or bet) is heads, and* ***m/(n+m)****, otherwise (when that choice is tails).*

Observe that **n = 1** and **m = 2** in the original Sleeping Beauty Problem.

*Proof*. By means of the same arguments as above, from Sections 4 and 5, when there shall be **n** awakenings with interview if the outcome of the coin flipping was heads, and **m** ones if that outcome was tails, we can conclude that the probability of getting heads [$P(H)$ = **1/2**] is *equally distributed* between all those **n** *heads-awakenings* [$P(H_1)$, $P(H_2)$, …, $P(H_n)$, where $P(H_i)$ is the probability of getting heads and awakening on day **i**, and all these $P(H_i)$ = **1/2n**, leading to $P(H) = P(H_1) + P(H_2) + \ldots + P(H_n)$ = **n/2n** = **1/2**], and that the probability of getting tails [$P(T)$ = 1/2, too] is also *equally distributed* between all the other **m** *tails-awakenings* [$P(T_1)$, $P(T_2)$, …, $P(T_m)$, where $P(T_i)$ is the probability of getting tails and awakening on day **i**, and all these $P(T_i)$ = **1/2m**, leading to $P(T) = P(T_1) + P(T_2) + \ldots + P(T_m)$ = **m/2m** = **1/2**], in *the SB's point of view* (for, when she awakes, without new information, she concludes that it can be any day from **1** to *Greater*(**n**, **m**) [if **n ≥ m**, then **n**, otherwise **m**], there is no way she knows which day she has wakened), where *the SB's point of view* is individualized in theorem above by means that amnesia-induced that makes her forget the previous awakenings, and then we can construct the two tables below, demonstrating and calculating those probabilities and ratios of expectations of the quantity of her hits:

(Observe that those distributions and collapses of SB's point of view probabilities in Section 4 [when SB does know which day she has wakened] appear also into the first table below, but only for additional clarifications on the probabilities involved in the problem, even though they are really not necessary in order to demonstrate the Theorem 6.1 above):

| **SB Does Not Know Which Day She Has Wakened** | |
|---|---|
| Whichever Awakening Day with Interview | $P(H_1) = P(H_2) = \ldots = \ldots = P(H_n) = \mathbf{1/2n}$ |
| | $P(T_1) = P(T_2) = \ldots = \ldots = \ldots = P(T_m) = \mathbf{1/2m}$ |
| | $P(H) = P(T) = \mathbf{1/2}$ |
| | $P(H) + P(T) = \mathbf{1}$ |
| **SB Does Know Which Day She Has Wakened** | |
| Day **i** Awakening Day (Where **i** ≤ *Less*(**n**, **m**)) | $P(H) = P(H_i) = \mathbf{1/2}$ and $P(H_u) = \mathbf{0}$ (for all **u ≠ i**) |
| | $P(T) = P(T_i) = \mathbf{1/2}$ and $P(T_u) = \mathbf{0}$ (for all **u ≠ i**) |
| | $P(H) = P(H_i) + P(H_u) = \mathbf{1/2} = P(T_i) + P(T_u) = P(T)$ (for all **u ≠ i**) |
| | $P(H) + P(T) = \mathbf{1}$ |
| Day **u** Awakening Day (Where **u** > *Less*(**n**, **m**)) When **n < m** | $P(H) = P(H_i) = \mathbf{0}$ (for all **i**) |
| | $P(T) = P(T_u) = \mathbf{1}$ and $P(T_i) = \mathbf{0}$ (for all **i ≠ u**) |
| | $P(T) = P(T_u) + P(T_i) = \mathbf{1}$ (for all **i ≠ u**) |
| | $P(H) + P(T) = \mathbf{1}$ |
| Day **u** Awakening Day (Where **u** > *Less*(**n**, **m**)) When **n > m** | $P(T) = P(T_i) = \mathbf{0}$ (for all **i**) |
| | $P(H) = P(H_u) = \mathbf{1}$ and $P(H_i) = \mathbf{0}$ (for all **i ≠ u**) |
| | $P(H) = P(H_u) + P(H_i) = \mathbf{1}$ (for all **i ≠ u**) |
| | $P(H) + P(T) = \mathbf{1}$ |

**Table 6.1** 1st Help for Visualizing the Demonstration of the Theorem 6.1 above

| Coin Landed | Probability (SB's Degree of Belief or Credence) | Qty of Awakenings w.Interview | SB's Choice (Guess or Bet) | Qty of Hits | Expectation of Qty of Hits | Ratio of Expectations |
|---|---|---|---|---|---|---|
| Heads | 1/2 | n | **Heads** | n | n*1/2+0*1/2 = n/2 | (n/2)/(n+m)/2 = **n/(n+m)** |
| Tails | 1/2 | m | | 0 | | |
| Heads | 1/2 | n | **Tails** | 0 | 0*1/2+m*1/2 = m/2 | (m/2)/(n+m)/2 = **m/(n+m)** |
| Tails | 1/2 | m | | m | | |

**Table 6.2** 2nd Help for Visualizing the Demonstration of the Theorem 6.1 above □

Notice that if, during some interview, SB was asked: "What is your credence (degree of belief) now for the proposition that today is the day **k**?" (where **k** ≤ *Greater*(**n**, **m**)), then, if **k** ≤ *Less*(**n**, **m**) [if **n < m**, then **n**, otherwise **m**], then she should answer **1/2n** + **1/2m** = (**n** + **m**)/**2nm** (resulting from: $P(H_k) + P(T_k)$), otherwise (if **k** > *Less*(**n**, **m**)), she should answer **1/(2***Greater*(**n**, **m**)) (resulting from: either $P(H_k)$ or $P(T_k)$, in case of either **n** or **m** is the greater one, respectively).

So, as *Greater*(**1**, **2**) = **2**, and *Less*(**1**, **2**) = **1**, if, during some interview in the original Sleeping Beauty Problem, SB was asked: "What is your credence (degree of belief) now for the proposition that today is Monday?", then she should answer (**1** + **2**)/(**2*1*2**) = **3/4**; and, if she was asked: "What is your credence (degree of belief) now for the proposition that today is Tuesday?", then she should answer **1**/(**2*2**) = **1/4**.

Interestingly, by the way, if, during some interview, it was told SB that the outcome of the coin flipping was heads, and then SB was asked: "What is your credence (degree of belief) now for the proposition that today is the day **k**?", then, if **k** ≤ **n**, then she should answer **1/n**, and if **k** > **n**, then she should answer **0**; on the other hand, if it was told SB that was tails, then, if **k** ≤ **m**, then she should answer **1/m**, and if **k** > **m**, then she should answer **0**.

So, if, during some interview in the original Sleeping Beauty Problem, it was told SB that the outcome of the coin flipping was heads, and then SB was asked: "What is your credence (degree of belief) now for the proposition that today is Monday?", then she should answer **1/1** = **1** (100% certainty), and if SB was asked: "What is your credence (degree of belief) now for the proposition that today is Tuesday?", then she should answer **0**; on the other hand, if it was told SB that was tails, and then SB was asked: "What is your credence (degree of belief) now for the proposition that today is Monday?", then she should answer **1/2**, and if SB was asked: "What is your credence (degree of belief) now for the proposition that today is Tuesday?", then she should answer **1/2**, too.

## 7. Confusion about Conditional Probabilities and the SB Problem

Several analysis on the Sleeping Beauty Problem, as in [3, 7], consider that the probabilities involved in it are conditional, but they are not, since there is no new information available to SB during the course of the experiment (whenever she awakes, she already knew that there would be at least one awakening, and she already knew also that there would be no way she knows which day she has wakened). Then conditional probabilities and Bayesian statistical methods [10] are not applicable herein, whatsoever.

That is, without new information at all, the degree of belief (credence) of SB about the outcome of that coin flipping, as not a mere arbitrary subjective concept, cannot change from Sunday to Monday and after that.

## 8. Conclusion

Consequently, it is expected that those two new versions of the Problem, that explanation about the origin of the thirder position confusion, the theorem above, and that clarification about the confusion on conditional probabilities and the SB Problem can help to find a possible consensual solution to the original Sleeping Beauty Problem, expanding our understanding about the self-locating belief in several areas of knowledge and science.

## 9. References

[1] From Wikipedia, the free encyclopedia, "Sleeping Beauty Problem", unpublished, available: https://en.wikipedia.org/wiki/Sleeping_Beauty_problem.

[2] A. Zuboff, *One Self: The Logic of Experience*. Inquiry: An Interdisciplinary Journal of Philosophy. 33 (1): 39–68, 1990.

[3] A. Elga, *Self-locating Belief and the Sleeping Beauty Problem*. Analysis. 60 (2): 143–147, 2000.

[4] D. Lewis, *Sleeping Beauty: reply to Elga*. Analysis. 61 (3): 171–76, 2001.

[5] From Wikipedia, the free encyclopedia, "Principle of Indifference", unpublished, available: https://en.wikipedia.org/wiki/Principle_of_indifference.

[6] E. T. Jaynes, *Probability Theory: The Logic of Science*. Cambridge University Press, 2003.

[7] From Terence Tao's blog *What's new*, "A Bayesian probability worksheet" [at Problem 2], unpublished, available: https://terrytao.wordpress.com/2022/10/07/a-bayesian-probability-worksheet/.

[8] From Commonplaces – a Journal of First-Year Writing at Davidson, "Heads or Tails: Always 50/50", available: https://commonplaces.davidson.edu/vol-3/heads-or-tails-always-5050/.

[9] O. Kallenberg, *Foundations of Modern Probability*. Springer-Verlag, 2nd ed., New York, 2002.

[10] Gelman, Andrew; Carlin, John B.; Stern, Hal S.; Dunson, David B.; Vehtari, Aki; Rubin, Donald B., *Bayesian Data Analysis*, 3rd ed., Chapman and Hall/CRC, 2013.